\documentclass[a4paper,10pt,oneside]{article}
\usepackage{amsmath}
\usepackage{amssymb}
\usepackage{calrsfs}
\usepackage{amstext}
\usepackage{epsfig}
\usepackage{color}
\usepackage[english]{babel}

\usepackage{esvect}
\usepackage{amsthm}
\usepackage{amscd}   
\usepackage{esvect}

\newcommand{\col}[2]{ \textcolor{#1}{#2}}

\newcommand{\dentro}[1]{$\mbox{Spivey\;\em \cite{spi}}.$}

\newcounter {eser}[section]

   \newtheorem{remark}{Remark}[section]
   
  \newtheorem{lemma}{Lemma}[section]
  \newtheorem{proposition}{Proposition}[section]
  \newtheorem{corollary}{Corollary}[section]
  \newtheorem{theorem}{Theorem}[section]
\newcommand{\N}{\mathbb{N}}

\begin{document}

\title{\bf Proofs of some theorems for binomial transform  and  Fibonacci powers} 
 \author{\bf Roberto S\'anchez-Peregrino,\\
Universit\`a  di Padova,\\ Dipartimento 
di Matematica: "Tullio Levi-Civita" ,\\ Via Trieste 63, I-35121 Padova, Italy.\\
E-mail: sanchez@math.unipd.it,peregrino84@gmail.com}
\date{\today}

\maketitle
\begin{abstract}  Our aim in writing  this paper is to answer to both V. E. Hoggatt, JR \cite{hogg} and  Wessner\cite{wess}  
 on the next question:  find $\sum_{k=0}^n\binom{n}{k}F_{[k]}^p$, for the case  $p\equiv 1\, mod\, 4$ and $p\equiv 3\, mod\, 4$. \par
The  case  $p\equiv 0\, mod \,4$ and $p\equiv 2\, mod\, 4$,  Wessner  has given an answer. In particular, we give another presentation, another proof of the paper of Wessner. 
Our method use, essentially, the paper of  Boyadzhiev\cite{boy} 
\end{abstract}
\vspace{20pt}
{\section{Introduction}}

Note that several papers about the evaluation  of $\sum_{k=0}^n(\pm 1)^k \binom{n}{k}F^p_{[k]}$ and similar sums have been written many 
years old, here  are  just some citations Carlitz et al\cite{quattro},  Wessner\cite{wess}.\par
The structure of this paper is as follows. In  Section 2 we give another representation  of the lemma and theorem of 
 Boyadzhiev \cite{boy2}. In section 3, we give a list of properties of Fibonacci and Lucas numbers.
We use it to give another representation of binomial sums of Fibonacci powers in the paper of Wessner
 \cite{wess}, i.e.
 $\sum_{k=0}^n(\pm 1)^k\binom{n}{k} F_{k}^{p}$ where $p\equiv 0(2) \pmod{4} $.
 We use these notations $F_{[n]}:=F_n$ and $L_{[n]}:=L_n$ 
for the Fibonacci 
numbers and  Lucas Numbers. In  section 4 we answer the question of V. E. Hoggatt, JR.\cite{hogg}\par
 	The authors hope that the above results will bring new life to this old problem. In a future paper we complete the cases, left by Hogartt:
 V. E. Hoggatt, Jr. \cite{hogg}
\par
{\bf Mathematics Subject Classification 2010:} Primary 11B65; 05A10.\par 
{\bf Keywords:} Binomial Transform, binomial identities, back difference operator,Fibonacci numbers and Lucas Numbers.
\newpage
{\section{Proofs}}
In this section, we give the other proof of paper \cite{boy}.
Given a sequence $\{a_n\}$, its binomial transform is the sequence $\{b_n\}$ defined by the formula 
\begin{eqnarray}
b_n=\sum_{k=0}^n\binom{n}{k}a_k,\label{eqn1}
\end{eqnarray}
with inversion
\begin{eqnarray}
a_n=\sum_{i=0}^n\binom{n}{k}(-1)^{n-k}b_k.\label{eqn2}
\end{eqnarray}
Writing for brevity $\nabla b_n=b_n-b_{n-1}$ and $\nabla ^0b_n=b_n$.
\begin{lemma}
\label{lemma1} Suppose the sequences $\{a_n\}$ and $\{b_n\}$ are defined by the formula (\ref{eqn1}).  
Then for every two  integers $0\leq m\leq n,$ we have

\begin{eqnarray}
\binom{n}{m}\nabla^m b_n=\sum_{j=0}^n\binom{n }{ j}\binom{j}{n-m}(-1)^{n-j}b_j\label{eqn3}
\end{eqnarray}
\text{or , in a shorter form}
\begin{eqnarray}
\nabla^m b_n=\sum_{k=0}^n\binom{m }{ k}(-1)^{k}b_{n-k}.\label{eqn4}
\end{eqnarray}
\end{lemma}
\proof 
The proof of  (\ref{eqn4}) is  by induction. \par 
When $m=0$ this is clearly true. We take $1\leq m< n$. We do the proof for $m+1.$
\begin{eqnarray*}
\nabla^{m+1} b_n=\sum_{k=0}^n\binom{m }{ k}(-1)^{k}\nabla b_{n-k}=\\
\sum_{k=1}^m\binom{m}{k}(-1)^kb_{n-k}+\binom{m}{0}b_n-\sum_{k=0}^{m-1}\binom{m}{k}(-1)^{k}b_{n-k-1}-(-1)^mb_{n-m-1}=\\
\sum_{k=1}^m\binom{m+1}{k}(-1)^{k}b_{n-k}+\binom{m}{0}b_n+(-1)^{m+1}b_{n-m-1}=
\sum_{k=0}^{m+1}\binom{m +1}{ k}(-1)^{k}b_{n-k}
\end{eqnarray*}
\endproof
\begin{lemma}
\label{lemma2} Suppose the sequences $\{a_n\}$ and $\{b_n\}$ are defined by the formula (\ref{eqn1}).Then for every two integers $0\leq m\leq n$ we have 
\begin{eqnarray}
\sum_{k=0}^n\binom{n}{k}\binom{k}{m}a_k=\binom{n }{ m}\nabla ^mb_n\label{eqn5}
\end{eqnarray}
\text{or , in a shorter form}
\begin{eqnarray}
\sum_{k=0}^n\binom{n-m}{k-m}a_k=\nabla ^mb_n.\label{eqn6}
\end{eqnarray}
\end{lemma}
\proof 
The RHS of  (\ref{eqn6}) is by (\ref{eqn4})  
\begin{eqnarray*}
\sum_{k=0}^n\binom{m }{ k}(-1)^{k}b_{n-k}=\sum_{k=0}^m\binom{m }{ k}(-1)^{k}b_{n-k}=\sum_{k=0}^m \sum_{l=0}^{n-k}\binom{m }{k}
\binom{n-k}{l}(-1)^{k}a_{l}=\\
\sum_{k=0}^m \sum_{l=0}^{n}\binom{m }{k}\binom{n-k}{l}(-1)^{k}a_{l}{=} 
\sum_{l=0}^{n}\sum_{k=0}^m \binom{m }{k}\binom{n-k}{l}(-1)^{k}a_{l}\stackrel{Gould:3.49}{=}\\
\sum_{l=0}^n\binom{n-m}{l-m}a_l.
\end{eqnarray*}
\endproof
\begin{lemma}
\label{lemma3} 
For any two integers $1\leq m\leq n.$ 
\begin{eqnarray}
\sum_{k=m}^n\binom{n}{k}\binom{k}{m}\frac{(-1)^k}{k}=\frac{(-1)^m}{m}.\label{eqn7}
\end{eqnarray}
\end{lemma}
\proof 

The LHS of  (\ref{eqn7}) is   
\begin{eqnarray*}
\sum_{k=m}^n\binom{n}{k}\binom{k}{m}\frac{(-1)^k}{k}=\binom{n}{m}\sum_{k=m}^n\binom{n-m}{k-m}\frac{(-1)^k}{k}=\binom{n}{m}(-1)^{m}\sum_{j=0}^{n-m}\binom{n-m}{j}\frac{(-1)^{j}}{m+j}=\\
=\binom{n}{m}\frac{(-1)^{m}}{m}\sum_{j=0}^{n-m}\binom{n-m}{j}\frac{(-1)^{j}m}{m+j}\stackrel{Gould:1.41}{=}
\binom{n}{m}\frac{(-1)^{m}}{m}\frac{1}{\binom{n}{n-m}}=\frac{(-1)^m}{m}.
\end{eqnarray*}
\endproof
\begin{theorem}
\label{theorem1}
Let $\{a_n\}$ and $\{c_n\}$ be two sequences and $b_n=\sum_{k=0}^n\binom{n}{k}a_k$ ,  $d_n=\sum_{k=0}^n\binom{n}{k}(-1)^{n-k}c_k$
then we have the identity
 \begin{equation}
\sum_{k=0}^n\binom{n}{k}a_kc_k=\sum_{m=0}^n\binom{n}{m}d_m\nabla^mb_n.\label{eqn8}\end{equation}
\text{or another form}
\begin{eqnarray}
\sum_{k=0}^n\binom{n}{k}a_kc_k=\sum_{k=0}^n(-1)^k\binom{n}{k}b_{n-k}\sum_{l=0}^n\binom{n-k}{l}d_{l+k}.\label{eqn81}
\end{eqnarray}
\end{theorem}
\proof 
The LHS of  (\ref{eqn8}) is 
\begin{eqnarray*}
\sum_{k=0}^n\binom{n}{k}a_kc_k=\sum_{k=0}^nd_n\left\{\sum_{k=0}^n\binom{n}{k}\binom{k}{m}a_k\right\}=\sum_{k=0}^n\binom{n}{m}d_n\left\{\sum_{k=0}^n\binom{n-m}{k-m}a_k\right\}=\\
\sum_{k=0}^n\binom{n}{m}d_n\nabla^mb_n.
\end{eqnarray*}
 Let $p(x)=\sum_{k=0}^n\binom{n}{k}e_kx^k$ then we define $f_n=p(1)=\sum_{k=0}^n\binom{n}{k}e_k.$
\begin{corollary}
\label{corollary1}
Let $\{a_n:=(-1x)^n\}$ and $\{c_n:=(-1)^ne_n\}$ be two sequences. We set   two sequences  $b_n=\sum_{k=0}^n(-1)^{n-k}\binom{n}{k}a_k$ , 
 $d_n=\sum_{k=0}^n\binom{n}{k}(-1)^{n-k}c_k.$ 

 \begin{equation}
\sum_{k=0}^n\binom{n}{k}e_kx^k=\sum_{n=0}^n\binom{n}{j}f_jx^j(1-x)^{n-j}=(1-x)^n\sum_{m=0}^n\binom{n}{j}f_j\left(\frac{x}{1-x}\right)^j.\label{eqn9}\end{equation}
\end{corollary}
\proof Then $d_n=(-1)^n\sum_{k=0}^n\binom{n}{k}e_k=(-1)^nf_n$.
 We have the identity  $b_n=(1-x)^n$. By (\ref{eqn4})
 we have $\Delta^mb_n=(1-x)^{n-m}(-x)^m$ and the LHS of equation(\ref{eqn9}) is the same of 
$\sum_{k=0}^n\binom{n}{k}(-1)^ka_k(-1)^ke_k=$ and   by theorem (\ref{theorem1})  is the same $\sum_{m=0}^n\binom{n}{m}d_m\Delta^m b_n =
(1-x)^n\sum_{m=0}^n\binom{n}{m}f_m\left(\frac{x}{1-x}\right)^m. $
\endproof
\begin{corollary}
\label{corollary2} Suppose the sequence $\{b_n\}$ is the binomial transform of the sequence $\{a_n\}.$ Then 

 \begin{equation}
\sum_{k=0}^n\binom{n}{k}a_kx^k=\sum_{m=0}^n\binom{n}{m}\nabla^mb_n(x-1)^{m}.\label{eqn91}\end{equation}
\end{corollary}

{\section{ Binomial transform of Fibonacci numbers}}
In this section we give another proof of paper \cite{wess}. Throughout this section, we denote the golden ratio, $\frac{1 + \sqrt{5}}{2}=x$ and write $y =\frac{ 1-\sqrt{5}}{2} = \frac{-1}{x}$,
\begin{remark} We have these relations between the golden ratio, $y$,  Fibonacci Numbers, and  Lucas numbers. 
Once this formul\ae  \, have been guessed the proofs are somewhat simpler.
\label{rem1}
\begin{eqnarray}
x^{4p}+1=x^{2p}L_{[2p]},\\
x^{4p}-1=x^{2p}\sqrt{5}F_{[2p]},\\
x^{4p-2}+1=x^{2p-1}\sqrt{5}F_{[2p-1]},\\
x^{4p+2}+1=x^{2p+1}\sqrt{5}F_{[2p+1]},\\
x^{4p-2}-1=x^{2p-1}L_{[2p-1]},\\
x^{4p+2}-1=x^{2p+1}L_{[2p+1]},\\
y^{4p}+1=y^{2p}L_{[2p]},\\
y^{4p}-1=-y^{2p}\sqrt{5}F_{[2p]},\\
y^{4p-2}+1=-y^{2p-1}\sqrt{5}F_{[2p-1]},\\
y^{4p+2}+1=-y^{2p+1}\sqrt{5}F_{[2p+1]},\\
y^{4p-2}-1=y^{2p-1}L_{[2p-1]},\\
y^{4p+2}-1=y^{2p+1}L_{[2p+1]}.
\end{eqnarray}
\end{remark}

\begin{proposition}
\label{propo11} 
Let $\{F_{[n]}\}$ be the  sequence  of Fibonacci Numbers then we have these relations 
 
\begin{eqnarray}
\sum_{k=0}^n\binom{n}{k}x^{kp}F_{[k]}=\frac{(x^{p+1}+1)^n-(-1)^n(x^{p-1}-1)^n}{\sqrt{5}}\label{eqn811}\\
\sum_{k=0}^n\binom{n}{k}(-x^p)^{k}F_{[k]}=\frac{(-1)^n(x^{p+1}-1)^n-(x^{p-1}+1)^n}{\sqrt{5}}\label{eqn812}\\
\sum_{k=0}^n\binom{n}{k}y^{kp}F_{[k]}=\frac{(-1)^n(y^{p-1}-1)^n-(y^{p+1}+1)^n}{\sqrt{5}}\label{eqn813}\\
\sum_{k=0}^n\binom{n}{k}(-y^p)^{k}F_{[k]}=\frac{(y^{p-1}+1)^n-(-1)^n(y^{p+1}-1)^n}{\sqrt{5}}\label{eqn814}
\end{eqnarray}
\end{proposition}
\proof  We give the proof of the equation (\ref{eqn811}),  the proofs the other equations use the likewise method.
 By the corollary(\ref{corollary2}) and the equation (\ref{eqn4}) we have the LHS of (\ref{eqn811})  is the same \par
\proof
\begin{eqnarray*}
\sum_{m=0}^n\binom{n}{m}\nabla^mF_{[n]}(x^p-1)^m= \sum_{m=0}^n\sum_{l=0}^{m}\binom{n}{m}\binom{m}{l}(-1)^lF_{[2n-2l]}(x^p-1)^m=\\
 \sum_{l=0}^n\binom{n}{l}(-1)^lF_{[2n-2l]}\sum_{m=0}^{m}\binom{n-l}{n-m}
(x^p-1)^m=
\sum_{l=0}^n\binom{n}{l}(-1)^lF_{[2n-2l]}(x^p)^n(\frac{x^p-1}{x^p})^l\\=\frac{x^{pn}}{\sqrt{5}}\sum_{l=0}^n(-1)^l(\frac{x^p-1}{x^p})^l(x^{2(n-l)}-y^{2(n-l)})=
\frac{(x^{p+1}+1)^n-(-1)^n(x^{p-1}-1)^n}{\sqrt{5}}.
\end{eqnarray*}
By the Proposition \ref{propo11}
and the Remark \ref{rem1}
\begin{eqnarray*}
\sum_{k=0}^n\binom{n}{k}F^{4p}_{[k]}=(\frac{1}{\sqrt{5}})^{4p-1}\left\{\sum_{r=0, r\; even}^{2p-1}\left[\binom{4p-1}{r}\frac{(x^{4p-2r}+1)^n-(-1)^n(x^{4p-2r-2}-1)^n}{\sqrt{5}}\right]+\right.\\
\sum_{r=0, r\; odd}^{2p-1}\left[\binom{4p-1}{r}\frac{(-1)^n(x^{4p-2r}+1)^n-(x^{4p-2r-2}-1)^n}{\sqrt{5}}\right]+\\
\sum_{r=0, r\; even}^{2p-1}\left[\binom{4p-1}{r}\frac{(y^{4p-2r}+1)^n-(-1)^n(y^{4p-2r-2}-1)^n}{\sqrt{5}}\right]+\\
\left.\sum_{r=0, r\; odd}^{2p-1}\left[\binom{4p-1}{r}\frac{(-1)^n(y^{4p-2r}+1)^n-(y^{4p-2r-2}-1)^n}{\sqrt{5}}\right]\right\}.
\end{eqnarray*}
\endproof
\begin{theorem}
\label{theorem2}
 {\bf {\mbox (Wessner)}}  We have

$\sum_{k=0}^n\binom{n}{k}F^{4p}_{[k]}=
\left\{\begin{array}{ll}
(\frac{1}{\sqrt{5}})^{4p}\left(
\binom{4p}{2p}2^n+\sum_{i=0}^{2p-1}(-1)^i\binom{4p}{i}L_{[2p-i]}^nL_{[(2p-i)n]}\right)&\mbox{if $n$ is even;}\\
(\frac{1}{\sqrt{5}})^{4p}\left(
\binom{4p}{2p}2^n+\sum_{i=0}^{2p-1}\binom{4p}{i}L_{[2p-i]}^nL_{[(2p-i)n]}\right)&\mbox{if $n$ is odd.}
\end{array}\right.$
\end{theorem}

\proof We have this relation:
\begin{eqnarray}\sum_{k=0}^n\binom{n}{k}F^{4p-1}_{[k]}F_{[k]}=
&&(\frac{1}{\sqrt{5}})^{4p-1}\left\{\sum_{k=0}^n\binom{n}{k}\sum_{r=0}^{2p-1}\binom{4p-1}{r}(-1)^r(x^{4p-1-2r}(-1)^r)^kF_{[k]}\right.\nonumber\\
&&+\left.\sum_{k=0}^n\binom{n}{k}\sum_{r=0}^{2p-1}\binom{4p-1}{r}(-1)^{r+1}(y^{4p-1-2r}(-1)^r)^kF_{[k]}\right\}=\nonumber\\
&&(\frac{1}{\sqrt{5}})^{4p-1}\left\{\sum_{k=0}^n\binom{n}{k}\sum_{r=0}^{p-1}\binom{4p-1}{2r}(x^{4p-1-4r})^kF_{[k]}\right.\nonumber\\
&&-(-1)^n\sum_{k=0}^n\binom{n}{k}\sum_{r=0}^{p-1}\binom{4p-1}{2r+1}((-x)^{4p-4r-3})^kF_{[k]}\nonumber\\
&&-\sum_{k=0}^n\binom{n}{k}\sum_{r=0}^{p-1}\binom{4p-1}{2r}(y^{4p-4r-1})^kF_{[k]}\nonumber\\
&&+(-1)^n\left.\sum_{k=0}^n\binom{n}{k}\sum_{r=0}\binom{4p-1}{2r+1}((-y)^{4p-4r-3})^kF_{[k]}\right\}=\nonumber\\
&&(\frac{1}{\sqrt{5}})^{4p-1}\left\{\sum_{r=0}^{p-1}\binom{4p-1}{2r}\frac{(x^{4p-4r}+1)^n-(-1)^n(x^{4p-4r-2}-1)^n}{\sqrt{5}}+\binom{4p-1}{2p-1}2^n\right.\nonumber\\
&&+\left.\sum_{r=0}^{p-1}\binom{4p-1}{2r+1}\frac{(y^{4p-4r}+1)^n-(-1)^n(y^{4p-4r-2}-1)^n}{\sqrt{5}}+\binom{4p-1}{2p-1}2^n\right\}\nonumber\\
&&=(\frac{1}{\sqrt{5}})^{4p}\left\{\sum_{r=0}^{p-1}\left(\binom{4p-1}{2r}{x^{2(p-r)n}L_{[2(p-r)]}^n-(-1)^n\binom{4p-1}{2r+1}x^{(2p-2(r+1))n}L_{[2p-(2r+1)]}^n}\right)\right.\nonumber\\
&&+\left.\sum_{r=0}^{p-1}\left(\binom{4p-1}{2r+1}y^{2(p-r)n}L_{[2(p-r)]}^n-(-1)^n\binom{4p-1}{2r+1}y^{(2p-(2r+1))n}L_{[(2p-(2r+1))n]}^n\right)+\binom{4p}{2p}2^n\right\}\nonumber\\ 
&&\stackrel{n\;even}{=}(\frac{1}{\sqrt{5}})^{4p}\left\{\sum_{r=0}^{2p-1}\binom{4p-1}{r}(-1)^rx^{2(p-r)n}L_{[2(p-r)]}^n+\sum_{r=0}^{2p-1}\binom{4p-1}{r}(-1)^ry^{2(p-r)n}L_{[2(p-r)]}^n\right.\nonumber\\
&&\left.+\binom{4p}{2p}2^n\right\}\nonumber\\
&&\stackrel{n\;odd}{=}(\frac{1}{\sqrt{5}})^{4p}\left\{\sum_{r=0}^{2p-1}\binom{4p-1}{r}x^{2(p-r)n}L_{[2(p-r)]}^n+\sum_{r=0}^{2p-1}\binom{4p-1}{r}y^{2(p-r)n}L_{[2(p-r)]}^n\right.\nonumber\\
&&\left.+\binom{4p}{2p}2^n\right\}\nonumber\\
\end{eqnarray}

\endproof
By the likewise methods, we can obtain the next theorems. We leave those for the interested reader.
\begin{theorem}
\label{theorem3}
 {\bf {\mbox (Wessner)}}  We have

$\sum_{k=0}^n\binom{n}{k}F^{4p+2}_{[k]}=
\left\{\begin{array}{ll}
(\frac{1}{\sqrt{5}})^{4p+2}\left(
\sum_{i=0}^{2p+1}(-1)^i\binom{4p+2}{i}(\sqrt{5}F_{[2p-i+1]})^nL_{[(2p-i+1)n]}\right)&\mbox{if $n$ is even;}\\
(\frac{1}{\sqrt{5}})^{4p+2}\left(
\sum_{i=0}^{2p+1}\binom{4p+2}{i}(\sqrt{5}F_{[2p-i+1]})^nL_{[(2p-i+1)n]}\right)&\mbox{if $n$ is odd.}
\end{array}\right.$
\end{theorem}
\begin{theorem}
\label{theorem4}
 {\bf {\mbox (Wessner)}}  We have

$\sum_{k=0}^n\binom{n}{k}(-1)^kF^{4p}_{[k]}=
\left\{\begin{array}{ll}
(\frac{1}{\sqrt{5}})^{4p}\left(
\sum_{i=0}^{2p}(-1)^i\binom{4p}{i}L_{[2p-i]n}(\sqrt{5}F_{[(2p-i)n]})^n\right)&\mbox{if $n$ is even;}\\
(\frac{1}{\sqrt{5}})^{4p-1}\left(
\sum_{i=0}^{2p}(-1)^i\binom{4p}{i}F_{[2p-i]n}(\sqrt{5}F_{[(2p-i)n]})^n\right)&\mbox{if $n$ is odd.}\\
\end{array}\right.$
\end{theorem}
\begin{theorem}
\label{theorem5}
 {\bf {\mbox (Wessner)}}  We have

$\sum_{k=0}^n\binom{n}{k}(-1)^kF^{4p+2}_{[k]}=
\left\{\begin{array}{ll}
(\frac{1}{\sqrt{5}})^{4p+2}\left(
\sum_{i=0}^{2p}(-1)^i\binom{4p+2}{i}(L_{[2p+1-i]})^nL_{[(2p+1-i)n]}+\binom{4p+2}{2p+1}2^n\right)&\mbox{if $n$ is even;}\\
-(\frac{1}{\sqrt{5}})^{4p+2}\left(
\sum_{i=0}^{2p}\binom{4p+2}{i}(L_{2p+1-i})^nL_{[(2p+1-i)n]}+\binom{4p+2}{2p+1}2^n\right)&\mbox{if $n$ is odd.}
\end{array}\right.$
\end{theorem}
{\section{Main results }}
In this section, we answer the question of Hoggatt and Wessner. We give the proof  the case $p\equiv 1 (mod\; 4)$\par 

\begin{lemma}
\label{lemma4}
 {\bf {\mbox (lemma4)}}
Let $n,c\in \N$ be. We set:  
$q(n,c)=\sum_{m=0}^n(-1)^m\binom{n+1}{2m+c+1},$
 then we have these relations  for $n\geq c\geq 0.$\par
$\begin{array}{ll}
q(n+1,c)=q(n,c-1)+q(n,c),\\
q(n+1,c)=\binom{n+1}{c}-q(n,c+1)+q(n,c).\\
\end{array}$
\end{lemma}
Using the lemma \ref{lemma4} yields a proof of next lemma,
\begin{lemma}
\label{lemma5}
 {\bf {\mbox (lemma5)}}
We set $g(n):=\sum_{j=0}^n(a+1)^{n-j}(b+1)^j$, with $ab=-1$ then we have this relation
$$g(n)=\sum_{c=1}^{n}q(n,c)(a^{j}+b^{j})+q(n,0).$$
\end{lemma}
\begin{lemma}
\label{lemma6}
 {\bf {\mbox (lemma6)}}
Let $n,c\in \N$ be. We set:  
$s(n,c)=(-1)^{n+c}\sum_{m=0}^n(-1)^m\binom{n+1}{2m+c+1}.$
 then we have these relations  for $n> c\geq 1.$\par
$\begin{array}{ll}
s(n+1,0)=s(n,1)-s(n,0)+(-1)^{n+1}\binom{n+1}{0},\\
s(n+1,c)=s(n,c-1)-s(n,c),\\
s(n+1,c)=(-1)^{n-c+1}\binom{n+1}{c}-s(n,c+1)-s(n,c).\\

\end{array}$
\end{lemma}
Using the lemma \ref{lemma6} yields a proof of next lemma,
\begin{lemma}
\label{lemma7}
 {\bf {\mbox (lemma7)}}
We set $g(n):=\sum_{j=0}^n(a-1)^{n-j}(b-1)^j$, with $ab=-1$ then we have this relation
$$g(n)=\sum_{c=1}^{n}s(n,c)(a^{j}+b^{j})+s(n,0).$$
\end{lemma}

\begin{theorem}
\label{theorem6}
 We have

$\sum_{k=0}^n\binom{n}{k}F^{4p+1}_{[k]}=
\left\{\begin{array}{ll}
(\frac{1}{\sqrt{5}})^{4p}\left(
\sum_{t=0}^{\col{red}{p-1}}\binom{4p+1}{2t}F_{[4p-4t+1]}\left\{\sum_{j=1}^{n-1}q(n-1,j)L_{[(4p-4t+1)j]}+q(n-1,0)\right\}-\right.\\
(-1)^n\left.
\sum_{t=0}^{p-1}\binom{4p+1}{2t+1}F_{[4p-4t-1]}\left\{\sum_{j=1}^{n-1}s(n-1,j)L_{[(4p-4t-1)j]}+s(n-1,0)\right\}+\binom{4p+1}{2p}F_{[2n]}\right).\\
\end{array}\right.$
\end{theorem}


\begin{theorem}
\label{theorem7}
We have

$\sum_{k=0}^n\binom{n}{k}F^{4p+3}_k=
\left\{\begin{array}{ll}
(\frac{1}{\sqrt{5}})^{4p+2}\left(
\sum_{t=0}^{\col{red}{p}}\binom{4p+3}{2t}F_{[4p-4t+3]}\left\{\sum_{j=1}^{n}q(n-1,j)L_{[(4p-4t+3)j]}+q(n-1,0)\right\}+\right.\\
-(-1)^n\left.
\sum_{t=0}^{p-1}\binom{4p+3}{2t+1}F_{[4p-4t+1]}\left\{\sum_{j=1}^{n-1}s(n-1,j)L_{[(4p-4t+1)j]}+s(n-1,0)\right\}+\binom{4p+3}{2p+1}F_{[n]}\right)\\
\end{array}\right.$
\end{theorem}


{\section{Acknowledgments }} I would like to thank K.N. Boyadzhiev for valuable discussions during the preparation of the paper, 
 for his advice and helpful comments

{ 
\end{document}